\documentclass[11pt]{article}
\usepackage[toc,page]{appendix}
\usepackage{mathrsfs}
\usepackage{amsfonts}
\usepackage{layout}
\usepackage[top=2cm, bottom=2cm, left=2.5cm, right=2.5cm]{geometry}
\usepackage{amssymb}
\usepackage{amsmath}
\usepackage{amsthm}
\usepackage[colorlinks=true, linkcolor=blue, citecolor=red]{hyperref}
\makeatletter

\@addtoreset{equation}{section}
\makeatother

\newtheorem{corollary}{Corollary}[section]

\newtheorem{proposition}{Proposition}[section]
\newtheorem{remark}{Remark}[section]

\newtheorem{theorem}{Theorem}[section]

\begin{document}

\title{New results on the associated Meixner, Charlier, Laguerre, and Krawtchouk polynomials}
\author{\textbf{Khalid Ahbli}}

\maketitle
\begin{center}
Faculty of Sciences of Agadir, Ibn Zohr University, Morocco\\
E-mail: ahbli.khalid@gmail.com
\vspace*{0.2mm}\vspace*{0.2mm}
\end{center}

\begin{abstract}We give new explicit representations as well as new generating functions for the associated Meixner, Charlier, Laguerre, and Krawtchouk polynomials. The obtained results are then used to derive new generating functions and convolution-type formulas of the corresponding classical polynomials. Some consequences of our results are also mentioned.\\

\end{abstract}
\section{Introduction}
A set of polynomials $\{P_n(x)\}$ is orthogonal if $P_n(x)$ is a polynomial of exact degree $n$ and there is a positive measure $d\mu$ on the real line with an infinite number of points increase and for which all the moments are finite such that $
\int_{\mathbb{R}}P_n(x)P_m(x)d\mu(x)=h_n \delta_{nm}, \ n,m=0,1,\cdots , \ (h_n>0).$ A necessary and sufficient condition for orthogonality \cite{Chihara} is that $\{P_n(x)\}$ satisfies the three-term recurrence
relation 
\begin{equation}\label{3TRR}
\begin{split}
& A_{n}P_{n+1}(x)=\left( B_{n}\,x+C_{n}\right)
P_{n}(x)-D_{n}P_{n-1}(x),\quad n=0,1,2,\cdots ,\\
& P_{-1}(x)=0,\text{ }P_{0}(x)=1,
\end{split}
\end{equation}
where $A_{n},\,B_{n},\,C_{n},\,D_{n}$ are real coefficients such that 
\begin{equation}
A_{n-1}B_{n-1}B_{n}D_{n}>0,\quad n\geq 1.  \label{coef_cond}
\end{equation}%
The associated orthogonal polynomials $\{P_{n}(x;\gamma )\}$ are defined by \eqref{3TRR}-\eqref{coef_cond} in terms of coefficients $A_{n+\gamma }$, $B_{n+\gamma }$%
, $C_{n+\gamma }$ and $D_{n+\gamma }$ where $\gamma \geq 0$ is a real association parameter. These polynomials have recently generated wide interest and some work has been devoted to the task of obtaining explicit representation (or closed-form expression) for these polynomials. This has been done successfully for the associated Jacobi polynomials \cite{Wimp1987} and their special cases, including the associated Gegenbauer \cite{Lew1993}, Laguerre and Hermite polynomials \cite{AW1984} and more recently for the associated Meixner-Pollaczek polynomials \cite{Luo2019}.  Associated orthogonal polynomials have many applications in diverse fields, such as queuing and inventory models, chemical kinetics, population dynamics, and quantum optics. For instance, it was shown in \cite{ILV_1988} that associated Meixner polynomials are birth and death process polynomials (a stationary Markov process) with rates $\lambda_n=c(n+\gamma+\beta)$, $\mu_n=n+\gamma$, $0<c<1$. In \cite{Ahbli2018,Ahbli2020}, associated Meixner-Pollaczek and Hermite polynomials were used to construct some new sets of nonlinear coherent states which play an important role in quantum optics. For other works on associated orthogonal polynomials and their applications see the articles \cite{Hend1990, IR1991, BD1967, BI1982, ILVW1990(a), ILV_1988, ILVW1990(b), Rahman1996, Rahman2000, GIM1991, Lew1995, Darke2009, Wimp1990,Ahbli2018,Ahbli2020}, the book \cite{IsmailBook} and their references. \\

In this paper, we will be concerned with the associated Meixner polynomials (AMPs), denoted $\mathscr{M}_{n}(x;\beta,c,\gamma )$, that correspond to coefficients $A_{n+\gamma }=c,\,B_{n+\gamma }=c-1,\,C_{n+\gamma }=(c+1)(n+\gamma )+\beta
c,\,D_{n+\gamma }=(n+\gamma )(n+\gamma +\beta -1)$.  These polynomials were first studied by Ismail et al. in \cite{ILV_1988} where a generating function was derived. They proved that the generating function of these polynomials has the following integral representation:
\begin{equation}\label{int_repr_gener_AMP}
\sum_{n=0}^{+\infty}\frac{(c\,t)^n}{(\gamma+1)_n}\mathscr{M}_n(x;\beta , c, \gamma )=\gamma\left(1-c\,t\right)^{-\beta-x}\left(1-t\right)^{x}\int_{0}^{1}u^{\gamma -1}\left(1-c\,t\,u\right)^{x+\beta-1}\left(1-t\,u\right)^{-x-1}d\,u ,
\end{equation}
valid for $\beta >0$, $\gamma \geq 0$ and $c\neq 1$ (see \cite[p.345(4.4)]{ILV_1988}).
We will show how these polynomials are connected to the associated Meixner-Pollaczek polynomials, $\mathscr{P}_n^{(\nu)}(x;\varphi, \gamma )$, obtained in \cite{Pollaczek1950}. Next, we use this connection relation (see Eq.\eqref{AMPP_AMP} below) and the explicit representation of the polynomials $\mathscr{P}_n^{(\nu)}(x;\varphi, \gamma )$, recently obtained by Luo and Raina in \cite{Luo2019}, to derive a new explicit expression for the polynomials $\mathscr{M}_{n}(x;\beta,c,\gamma )$. The resulting formula will then be used to obtain a new generating function of these polynomials. Furthermore, we exploit the obtained results to derive new explicit representations as well as generating functions for the associated Charlier, Laguerre, and Krawtchouk polynomials. Different interesting identities, generating functions and convolution-type formulas (in terms of classical polynomials) are also established. Some proofs of the stated results are postponed to Sect. \ref{proofs}. A conclusion is also given at the end of this paper. 

\section{A new explicit formulae and generating function for the AMPs}\label{Sect_AMP}
%\subsection{A closed formula in terms of ${}_4F_3(1)$}
In this section, we will give a new explicit representation for the AMPs in terms of terminating ${}_4F_3(1)$-series which we deduced from the one of the associated Meixner-Pollaczek polynomials obtained in \cite{Luo2019}. The derived representation is then used to compute a new generating function of these polynomials. 
\subsection{Explicit representation of the AMPs}
The AMPs $\mathscr{M}_{n}(x;\beta
,c,\gamma )$ satisfy the three-term recurrence relation:
\begin{equation}\label{AM_RR}
\begin{split}
c\mathscr{M}_{n+1}(x;\beta ,c,\gamma )=[(c-1)x+(c+1)(n+\gamma )& +\beta c]%
\mathscr{M}_{n}(x;\beta ,c,\gamma ) \\
& -(n+\gamma )(n+\gamma +\beta -1)\mathscr{M}_{n-1}(x;\beta ,c,\gamma ),
\end{split}
\end{equation}
with initial conditions $\mathscr{M}_{-1}(x;\beta ,c,\gamma )=0$ and $\mathscr{M}_{0}(x;\beta ,c,\gamma )=1$. The conditions for orthogonality \eqref{coef_cond} require that $c>0,\,c\neq 1$ and $\gamma +\beta >0$. Now, we proceed to establish an explicit representation for the AMPs.
Our main tools here are the explicit representation for associated Meixner-Pollaczek polynomials $\mathscr{P}_n^{(\nu)}(x;\varphi, \gamma )$ obtained Luo and Raina in \cite[p.3(1.8)]{Luo2019},
\begin{eqnarray}\label{AMPP_4F3}
\begin{split}
\mathscr{P}_n^{(\nu)}(x;\varphi, \gamma )=\frac{(2\nu+\gamma)_n}{n!}\sum_{k=0}^{n}e^{i(n-k)\varphi}\left(e^{i\varphi}-e^{-i\varphi}\right)^k\frac{(-n)_k(\gamma+\nu+ix)_k}{(\gamma +1 )_k(\gamma+2\nu)_k}\\ 
\times {}_4F_3\left( 
\begin{array}{c}
k-n,\,\gamma+\nu+ix+k ,\, \gamma+2\nu -1,\, \gamma\ \\ 
 \gamma+\nu +ix,\, \gamma+2\nu +k, \, \gamma+1+k
\end{array}
\big| 1\right),
\end{split}
\end{eqnarray} 
 and the following connection relation:
\begin{equation}
\mathscr{P}_n^{(\nu)}(x;\varphi, \gamma )=\frac{e^{-in\varphi}}{(\gamma +1)_n}\mathscr{M}_n(ix-\nu ;2\nu,e^{-2i\varphi } ,\gamma ),\label{AMPP_AMP}
\end{equation}
which can be checked by comparing the recurrence relation \eqref{AM_RR} of $\mathscr{M}_n(x;\beta , c, \gamma )$ with the recurrence relation of $\mathscr{P}_n^{(\nu)}(x;\varphi, \gamma )$ given in \cite[p.2256]{Pollaczek1950} by
\begin{eqnarray}
\begin{split}\label{AMP_RR}
(n+\gamma+1)\mathscr{P}_{n+1}^{(\nu)}(x;\varphi, \gamma ) = 2[(n+\gamma+\nu)\cos\varphi +& x\sin\varphi ] \mathscr{P}_{n}^{(\nu)}(x;\varphi, \gamma ) \\
&-(n + \gamma+2\nu-1)\mathscr{P}_{n-1}^{(\nu)}(x;\varphi, \gamma ).
\end{split}
\end{eqnarray}
Here, $(a)_{k}$ is the Pochhammer symbol defined by $(a)_{0}=1$, $(a)_{k}=a(a+1)\cdots (a+k-1)$, $k\geq 1$. We are now in position to state the following result. 
\begin{proposition}\label{thm_CF_AMP}
We have the following explicit representation for the AMPs 
\begin{eqnarray}\label{AMP_4F3}
\begin{split}
\mathscr{M}_n(x;\beta , c, \gamma )=c^{-n}\frac{(\gamma +1 )_n(\gamma+\beta)_n}{n!}\sum_{k=0}^{n}(1-c)^k\frac{(-n)_k(\gamma+\beta+x)_k}{(\gamma +1 )_k(\gamma+\beta)_k}\\ 
\times {}_4F_3\left( 
\begin{array}{c}
k-n,\,\gamma+\beta+x+k ,\, \gamma+\beta -1,\, \gamma\ \\ 
 \gamma+\beta +x,\, \gamma+\beta +k, \, \gamma+1+k
\end{array}
\big| 1\right),
\end{split}
\end{eqnarray} 
valid for $c>0,\,c\neq 1$ and $\gamma +\beta >0$. 
\end{proposition}
In view of the relation $\mathscr{M}_{n}(x;\beta, c,\gamma)=c^{-n}\mathscr{M}_{n}(-\beta-x;\beta,
c^{-1},\gamma)$, which follows directly from the recurrence relation, we obtain from \eqref{AMP_4F3} another expression of the AMPs given by
\begin{eqnarray}\label{AMP_4F3_2}
\begin{split}
\mathscr{M}_n(x;\beta , c, \gamma )=\frac{(\gamma +1 )_n(\gamma+\beta)_n}{n!}\sum_{k=0}^{n}\tilde{c}^k\frac{(-n)_k(\gamma-x)_k}{(\gamma +1 )_k(\gamma+\beta)_k} {}_4F_3\left( 
\begin{array}{c}
k-n,\,\gamma-x+k ,\, \gamma+\beta -1,\, \gamma\ \\ 
 \gamma-x,\, \gamma+\beta +k, \, \gamma+1+k
\end{array}
\big| 1\right)
\end{split}
\end{eqnarray}
where $\tilde{c}=\frac{c-1}{c}$. The case $\gamma =0$ in \eqref{AMP_4F3_2} corresponds to the Meixner polynomials \cite[p.176 (3.5)]{Chihara}:
\begin{equation}\label{Meixner_poly}
M_{n}(x;\beta, c)=(\beta)_n\, {}_2F{}_1\left( 
\begin{array}{c}
-n,\, -x \\ 
\beta
\end{array}
\left|1-\frac{1}{c}\right.\right).
\end{equation}
In the special case $\gamma+\beta=1$, the polynomials $\mathscr{M}_n(x;\beta , c, \gamma )$ reduce to a particular case of Meixner polynomials. In fact, when $\gamma+\beta=1$ Eq.\eqref{AMP_4F3_2} reduces to 
\begin{eqnarray}\label{MP_2}
\mathscr{M}_n(x;\beta , c, 1-\beta )=(2-\beta)_n \sum_{k=0}^{n}\frac{(-n)_k(1-\beta-x)_k}{(2-\beta )_k}\frac{\tilde{c}^k}{k!}=(2-\beta)_n \, {}_2F{}_1\left( 
\begin{array}{c}
-n,\, 1-\beta-x \\ 
2-\beta
\end{array}
\left|\tilde{c}\right.\right),
\end{eqnarray}
and hence
\begin{equation}
\mathscr{M}_n(x;\beta , c, 1-\beta )=M_{n}(x+\beta-1;2-\beta, c).
\end{equation}
The above connection formula is known; see \cite[Eq.(11)]{LRV_1996}.
\begin{corollary}\label{corollary_sum_4F3}
The following finite sum
\begin{eqnarray}  \label{sum_fin_4F3}
\begin{split}
&\sum_{k=0}^{n}t^k \frac{(-n)_k(a+y)_k}{(a +1 )_k(b+1)_k}{}_4F_3\left( 
\begin{array}{c}
k-n,\,a+y+k ,\, a,\, b \  \\ 
a+y,\, b+1 +k, \, a+1+k%
\end{array}
\big| 1\right)& \\
&=\frac{n!}{b-a}\left\{ \frac{b}{(a +1 )_n} \, {}_2F_1\left( 
\begin{array}{c}
1-y,a \\ 
a-b+1%
\end{array}
\big|t\right){}_2F_1\left( 
\begin{array}{c}
y,-n-a \\ 
b-a+1%
\end{array}
\big|t\right)\right. & \\
&\left. -(1-t)^{n+1}\frac{a}{(b+1)_n}\,{}_2F_1\left( 
\begin{array}{c}
1-y,n+a +1 \\ 
a-b+1%
\end{array}
\big|t\right){}_2F_1\left( 
\begin{array}{c}
y,1-a \\ 
b-a+1%
\end{array}
\big|t\right)\right\} &
\end{split}%
\end{eqnarray}
is valid for $a>0,\, b>-1,\, b\neq 0$ and $b-a\neq 0,1,2,...$ .
\end{corollary}
\noindent For special values of parameters ($y=1$ or $t=0$), we obtain: 
\begin{equation*}
\begin{split}
\sum_{k=0}^{n}t^{k}\frac{(-n)_{k}}{(b+1)_{k}}{}_{3}F_{2}\left( 
\begin{array}{c}
k-n,\,a,\,b\  \\ 
a+1,\,b+1+k%
\end{array}%
\big|1\right) =& \frac{n!}{b-a}\left\{ \frac{b}{(a+1)_{n}}%
\,{}_{2}F_{1}\left( 
\begin{array}{c}
1,-n-a \\ 
b-a+1%
\end{array}%
\big|t\right) \right. \\
& \left. -(1-t)^{n+1}\frac{a}{(b+1)_{n}}\,{}_{2}F_{1}\left( 
\begin{array}{c}
1,1-a \\ 
b-a+1%
\end{array}%
\big|t\right) \right\} ,
\end{split}%
\end{equation*}%
and
\begin{equation}
{}_{3}F_{2}\left( 
\begin{array}{c}
-n,\,a,\,b\  \\ 
a+1,\,b+1%
\end{array}%
\big|1\right) =\frac{n!}{b-a}\left\{ \frac{b}{(a+1)_{n}}-\frac{a}{(b+1)_{n}}%
\right\}.  \label{3F2_pocha}
\end{equation}
%\textcolor{red}{
%A generalization of \eqref{3F2_pocha}, for $m\in\mathbb{N}$, was proved in \cite[p.37(5.7)]{MP2012} and is given by 
%\begin{equation*}
%{}_{3}F_{2}\left( 
%\begin{array}{c}
%-s,\,a,\,b\  \\ 
%a+m,\,b+1%
%\end{array}%
%\big|1\right) =\frac{(a)_{m}}{(a-b)_{m}}\left\{ \frac{(1)_{s}}{(1+b)_{s}}-%
%\frac{b}{a}\sum\limits_{l=0}^{m-1}\frac{(a-b)_{l}(1+l)_{s}}{%
%(1+a)_{l}(1+a+l)_{s}}\right\} .
%\end{equation*}
%}

%\begin{remark}
%If we set 
%\begin{equation}
%\textbf{M}_n^{\beta , \gamma}(x)=\lim\limits_{c\to +\infty}\mathscr{M}_n(x;\beta , c, \gamma )
%\end{equation}
%\end{remark}
\begin{remark}
Note that the only known expression of these polynomials was obtained in a quadratic form (cross products) made out of Gauss hypergeometric functions  ${}_2F_1$ as (see \cite[Eq(21)]{LRV_1996}):
\begin{eqnarray}\label{m_n}
\begin{split}
\mathscr{M}_n(x;\beta ,c,\gamma)=\frac{1}{\beta -1}\left\{ (\gamma +\beta -1)_{n+1} \, {}_2F_1\left( 
\begin{array}{c}
x+1,\gamma \\ 
2-\beta
\end{array}
\big|\tilde{c}\right){}_2F_1\left( 
\begin{array}{c}
-x,-n-\gamma \\ 
\beta
\end{array}
\big|\tilde{c}\right)\right.\\
\left. -(\gamma)_{n+1}\,{}_2F_1\left( 
\begin{array}{c}
x+\beta,\gamma+\beta -1 \\ 
\beta
\end{array}
\big|\tilde{c}\right){}_2F_1\left( 
\begin{array}{c}
1-\beta-x,-n-\gamma-\beta +1 \\ 
2-\beta
\end{array}
\big|\tilde{c}\right)\right\}
\end{split}
\end{eqnarray}
where $\tilde{c}=\frac{c-1}{c}$. This
relation is valid under restrictions $\gamma >0$, $\gamma +\beta \in \mathbb{R}_+^*\backslash \{1\}$ and $\beta \notin \mathbb{N}^*$.\\
\end{remark}

\subsection{A generating function for the AMPs}
Let us first notice that the generating function \eqref{int_repr_gener_AMP} may be deduced from the generating function for the associated Meixner-Pollaczek polynomials derived in \cite[p.3(1.7)]{Luo2019} with the help of the connection relation \eqref{AMPP_AMP}. Precisely, we have
\begin{eqnarray}\label{gen_AMP}
\sum_{n=0}^{+\infty}\frac{(c\,t)^n}{(\gamma+1)_n}\mathscr{M}_n(x;\beta , c, \gamma )=\left(1-c\,t\right)^{-\beta-x}\left(1-t\right)^{x}F_1\left[\gamma , 1-\beta-x, 1+x ; \gamma +1 ; c\,t, t \right],
\end{eqnarray} 
where $F_1$ is the first Appell hypergeometric function defined by \eqref{Appel_F1}. Then, the integral in the R.H.S of \eqref{int_repr_gener_AMP} is just the integral representation of the Appell function $F_1$ (see \cite[p.77(4)]{Bailey1964}). 
%\textcolor{red}{Letting $c=1$ in \eqref{gen_AMP} and applying the following reduction formula for $F_1$ \cite[p.79(1)]{Bailey1964}:
%\begin{eqnarray}
%F_1\left[\alpha_1 ; \,\lambda_1 ,\, \lambda_2 ;\,\alpha_2 ;\, t,\, t \right]={}_2F_1\left( 
%\begin{array}{c}
%\alpha_1 ,\, \lambda_1 +\lambda_2\ \\ 
% \alpha_2
%\end{array}
%\big| t\right),
%\end{eqnarray}
%with $\alpha_1=\gamma,\, \lambda_1=1-\beta-x,\, \lambda_2=1+x$ and $\alpha_2=\gamma+1$, gives 
%\begin{eqnarray}
%\sum_{n=0}^{+\infty}\frac{(\gamma+\beta)_n }{n!} {}_3F_2\left( 
%\begin{array}{c}
%-n ,\, \gamma+\beta-1,\, \gamma\ \\ 
% \gamma+\beta ,\, \gamma+1
%\end{array}
%\big| 1\right)\,t^n=(1-t)^{-\beta}{}_2F_1\left( 
%\begin{array}{c}
%2-\beta,\, \gamma\ \\ 
% \gamma+1
%\end{array}
%\big| t\right),\quad |t|<1,
%\end{eqnarray}
%which can be identified as a special case of the formula given in \cite[p.107(14)]{Sr-Ma1984}.}
In addition, if we put $\gamma=0$ in \eqref{gen_AMP}, we recover the following generating function of Meixner polynomials
\begin{eqnarray}\label{gener_MPP}
\sum_{n=0}^{+\infty}\frac{(c\,t)^n}{n!}M_n(x;\beta , c)=\left(1-c\,t\right)^{-\beta-x}\left(1-t\right)^{x}.
\end{eqnarray}   
Substitute \eqref{gener_MPP} into \eqref{int_repr_gener_AMP} with $(t,\beta,x)=(tu,2-\beta, -x-1)$, integrate the obtained expression, and then use the Rainville formula \cite[p.101]{Sr-Ma1984},
\begin{equation}\label{RF}
\sum_{n=0}^\infty\sum_{k=0}^\infty A(k,n)=\sum_{n=0}^\infty\sum_{k=0}^n A(k,n-k),
\end{equation}
to get an expression of $\mathscr{M}_n(x;\beta , c, \gamma )$ in terms of $M_{n}(x;\beta , c)$ as follows .
\begin{proposition} Let $\gamma >0$. Then, we have the convolution-type formula 
\begin{equation}
\frac{n!\, \mathscr{M}_n(x;\beta , c, \gamma )}{(\gamma+1)_n}=\sum_{k=0}^{n} \displaystyle \binom{n}{k} \frac{\gamma }{(k+\gamma)}M_{n-k}(x;\beta , c)M_{k}(-x-1;2-\beta , c).
\end{equation} 
\noindent Moreover, we have the generating function
\begin{equation}\label{gen_MP}
\sum_{n=0}^{+\infty} \frac{\gamma \,(c\,t)^n}{(n+\gamma)n!}M_{n}(x;\beta , c)=F_1\left[\gamma , x+\beta , -x ; \gamma +1 ; c\,t, t \right],\quad |t|< |c|^{-1},\, (\gamma \text{ arbitrary}).
\end{equation}
\end{proposition}
The last formula is obtained by direct computation with the help of generating functions \eqref{gen_AMP} and \eqref{gener_MPP}.\\

We end this subsection by stating a new generating function for the AMPs (see Sect. \ref{proofs} for a proof).
\begin{theorem}\label{thm_GF_AMP}
A generating function for the AMPs is given by 
\begin{eqnarray}\label{gen2_AMP}
\sum_{n=0}^{+\infty}\frac{(c\,t)^n}{(\gamma+\beta)_n}\mathscr{M}_n(x;\beta , c, \gamma )=\left(1-c\,t\right)^{-1}F_1\left[1 , \gamma , -x ; \gamma +\beta ; t, \frac{t(1-c)}{1-ct} \right]
\end{eqnarray} 
for $|t|<\min \{1,\, c^{-1}\}$, in terms of Appell function $F_1$ defined by \eqref{Appel_F1}. In particular, we have
\begin{eqnarray}\label{gener2_MPP}
\sum_{n=0}^{+\infty}M_n(x;\beta , c)t^n=\left(1-t\right)^{-1}{}_2F_1\left( 
\begin{array}{c}
1,\, -x \ \\ 
\beta
\end{array}
\left| \frac{t(1-c)}{c(1-t)}\right.\right).
\end{eqnarray} 
\end{theorem}

\section{The associated Charlier polynomials}\label{sect_ACP}
The associated Charlier polynomials (ACPs) are defined by the three-term recurrence relation:
\begin{eqnarray}\label{ACP_RR}
a\,\mathscr{C}_{n+1}(x;a,\gamma )=(n+\gamma +a -x)\mathscr{C}_{n}(x;a,\gamma )-(n+\gamma)\mathscr{C}_{n-1}(x;a,\gamma )
\end{eqnarray}
with initial conditions $\mathscr{C}_{-1}(x;a,\gamma )=0$ and $\mathscr{C}_{0}(x;a,\gamma )=1$. Applying the criterion \eqref{coef_cond} we find that $\mathscr{C}_{n}(x;a,\gamma )$ are orthogonal if $a>0$.
By comparing the recurrence relations of the AMPs and ACPs, it is easy to check the limiting relation (see \cite[(35)]{LRV_1996}): 
\begin{equation}
\mathscr{C}_{n}(x;a,\gamma )=\lim\limits_{\beta \rightarrow \infty }\frac{1}{%
(\gamma +\beta )_{n}}\mathscr{M}_{n}\left( x;\beta ,\frac{a}{a+\beta }%
,\gamma \right).  \label{LR_CM}
\end{equation}%
Representation \eqref{AMP_4F3_2} provides for ACPs two interesting explicit formulas. The first one is obtained thanks to the relation \eqref{LR_CM}, upon using the limit 
\begin{equation}  \label{limit_pocha}
\lim\limits_{a\to \infty}\frac{(ax)^k}{(ay+b)_k}=\left(\frac{x}{y}\right)^k,
\ k\geq 1,
\end{equation}
where $x,y,b$ are fixed, and is given by
\begin{eqnarray}  \label{ACP_exp}
\mathscr{C}_n(x;a,\gamma )=\frac{(\gamma +1 )_n}{n!}\sum_{k=0}^{n}(-a)^{-k}%
\frac{(-n)_k(\gamma-x)_k}{(\gamma +1 )_k} {}_3F_2\left( 
\begin{array}{c}
k-n,\,\gamma-x+k ,\, \gamma\  \\ 
\gamma-x,\, \gamma+k+1%
\end{array}
\big| 1\right),\  a> 0.
\end{eqnarray}
The second formula is obtained from \eqref{ACP_exp} by applying the transformation (\cite[p.142]{AAR1999}):
\begin{eqnarray}  \label{3F2_trans}
{}_3F{}_2\left( 
\begin{array}{c}
-m,\,a,\,b \\ 
c,\,d%
\end{array}
\big|1\right)=\frac{(c-a)_m}{(c)_m} {}_3F_2\left( 
\begin{array}{c}
-m,\,a,\,d-b \\ 
a-c+1-m,\, d%
\end{array}
\big|1\right).
\end{eqnarray}
It is given by
\begin{eqnarray}  \label{ACP_exp2}
\mathscr{C}_n(x;a,\gamma )=\sum_{k=0}^{n}(-a)^{-k}\frac{(-n)_k(\gamma-x)_k }{k!}
\ {}_3F_2\left( 
\begin{array}{c}
-k,\,\gamma ,\, k-n \  \\ 
-n,\, \gamma-x%
\end{array}
\big| 1\right).
\end{eqnarray}
%\subsection{A generating function for associated Charlier polynomials}
Following the same lines as for the AMPs, we derive the generating function for the ACPs. Precisely, we have the following result (see Sect. \ref{proofs} for a proof).
\begin{proposition}\label{prop_GF_ACP}
A generating function for the ACPs is given by
\begin{eqnarray}\label{gener_ACP}
\sum_{n=0}^{+\infty}\frac{t^n}{(\gamma+1)_n}\mathscr{C}_n(x;a, \gamma )=e^{t}\left(1-\frac{t}{a}\right)^{x}\Phi_1\left[\gamma , x+1 ; \gamma +1 ; \frac{t}{a}, -t \right], \quad |t|<|a|,
\end{eqnarray} 
where $\Phi_1$ is the Humbert's confluent hypergeometric function defined by \eqref{Phi_1F1}.\\
In particular, we have
\begin{eqnarray}\label{gener_CP}
\sum_{n=0}^{+\infty}\frac{t^n}{n!}C_n(x;a)=e^{t}\left(1-\frac{t}{a}\right)^{x}.
\end{eqnarray} 
\end{proposition}

Let us now apply an alternative method to compute the generating function of $\mathscr{C}_n(x;a, \gamma )$. Multiplying both the sides of \eqref{ACP_RR} by $t^n/(\gamma+1)_n$ and then summing with respect to $n$ from $0$ to
$+\infty$, to get the following first-order differential equation
\begin{equation}\label{equ-diff-ACP}
t(a-t)\frac{\partial}{\partial t}\mathscr{G}_{\gamma}(x,t)+\left[t^2+(x-a-\gamma )t+a\gamma \right]\mathscr{G}_{\gamma}(x,t)=a\gamma 
\end{equation}
where $\mathscr{G}_{\gamma}(x,t)$ denotes the L.H.S of \eqref{gener_ACP}. 
Next, we make the substitution $\mathscr{G}_{\gamma}(x,t)=e^{t}\left(1-t/a\right)^{x}\mathscr{H}_{\gamma}(x,t)$
and obtain for $\mathscr{H}_{\gamma}(x,t)$ the following differential equation:
\begin{equation}
t\frac{\partial}{\partial t}\mathscr{H}_{\gamma}(x,t)+\gamma \mathscr{H}_{\gamma}(x,t)=\gamma e^{-t}\left(1-\frac{t}{a}\right)^{-x-1}.
\end{equation}
The solution of the above equation, after taking into account the condition $\mathscr{G}_{\gamma}(x,0)=1$, is 
\begin{equation}
\mathscr{H}_{\gamma}(x,t)=\gamma \int_{0}^{1}u^{\gamma-1}e^{-u\,t}\left(1-\frac{u\,t}{a}\right)^{-x-1}du, \quad |t|<\min \{1,\, |a|\}.
\end{equation}
Thus, $\mathscr{G}_{\gamma}(x,t)$ has the following integral representation 
\begin{equation}\label{int_repr_GF_ACP}
\mathscr{G}_{\gamma}(x,t)=\gamma  e^{t}\left(1-\frac{t}{a}\right)^{x}\int_{0}^{1}u^{\gamma-1}e^{-t\,u}\left(1-\frac{t\,u}{a}\right)^{-x-1}du.
\end{equation}
%\textcolor{red}{Note here that in order to simplify the differential equation \eqref{equ-diff-ACP}, we have made a substitution of type $\mathscr{G}_{\gamma}(x,t)=G(x,t)\mathscr{H}_{\gamma}(x,t)$ which is nothing than writing the generating function of the associated polynomials as a product of the (analogous) generating function of (corresponding) classical polynomials and another function which is now easy to find and satisfies $\mathscr{H}_{0}(x,t)=1$. We believe that this is the best substitution we may perform when we look for generating functions of associated orthogonal polynomials. } \\
Again, substitute \eqref{gener_CP} into \eqref{int_repr_GF_ACP} with $(t,\beta,x)=(-tu,-a, -x-1)$, integrate the obtained expression, and then use \eqref{RF}, to get the following corollary.
\begin{corollary} Let $\gamma >0$. Then we have
\begin{equation}
\frac{n!\,\mathscr{C}_n(x;a, \gamma )}{(\gamma+1)_n}=\sum\limits_{k=0}^{n}\displaystyle \binom{n}{k}\frac{\gamma(-1)^k}{(\gamma +k)}C_{n-k}(x;a)C_k(-x-1;-a).
\end{equation}
Consequently, we get the following generating function
\begin{eqnarray}\label{gener_2_CP}
\sum_{n=0}^{+\infty}\frac{\gamma}{\gamma+n}C_n(x;a)\frac{t^n}{n!}=\Phi_1\left[\gamma , -x ; \gamma +1 ; \frac{t}{a}, t \right], \quad |t|<|a|, \, (\gamma \text{ arbitrary}).
\end{eqnarray} 
\end{corollary}
\section{The associated Laguerre polynomials}\label{sect_ALP}
The associated Laguerre polynomials (ALPs) can be defined by the three-term recurrence relation:
\begin{eqnarray}
\begin{split}\label{AL_RR}
(n+\gamma+1)\mathscr{L}_{n+1}^{(\alpha)}(x;\gamma ) = [2(n+\gamma)+\alpha+1-x ] \mathscr{L}_n^{(\alpha)}(x;\gamma )-(n + \gamma+\alpha)\mathscr{L}_{n-1}^{(\alpha)}(x;\gamma )
\end{split}
\end{eqnarray}
with initial conditions $\mathscr{L}_{-1}^{(\alpha)}(x;\gamma )=0$ and $\mathscr{L}_0^{(\alpha)}(x;\gamma )=1$. Again, with the help of the criterion \eqref{coef_cond}, we show that $\mathscr{L}_n^{(\alpha)}(x;\gamma )$ are orthogonal if and only if $\alpha +\gamma >-1$.
By comparing \eqref{AM_RR} with \eqref{AL_RR}  it is easy to check the limit relation (see \cite[(39)]{LRV_1996}): 
\begin{equation}
\mathscr{L}_{n}^{(\alpha )}(x;\gamma )=\lim\limits_{c\rightarrow 1}\frac{1}{%
(\gamma +1)_{n}}\mathscr{M}_{n}\left( \frac{x}{1-c};\alpha +1,c,\gamma
\right) .  \label{LR_ALP-AMP}
\end{equation}%
This relation allows us to write an explicit formula for the ALPs as follows
\begin{equation}  \label{ALP_exp}
\mathscr{L}_n^{(\alpha)}(x;\gamma )=\frac{(\gamma +\alpha+1 )_n}{n!}%
\sum_{k=0}^{n}\frac{(-n)_k \, x^k}{(\gamma +1 )_k(\gamma +\alpha+1 )_k} {}%
_3F_2\left( 
\begin{array}{c}
k-n,\,\gamma+\alpha ,\, \gamma\  \\ 
\gamma+\alpha+k+1,\, \gamma+k+1%
\end{array}
\big| 1\right)
\end{equation}
where $\alpha >-1$.
Using the transformation %
\eqref{3F2_trans}, these polynomials can be rewritten as
follows (see also \cite[Eq(1.34)]{Rahman2000}): 
\begin{equation}
\mathscr{L}_n^{(\alpha)}(x;\gamma )=\frac{(\alpha+1 )_n}{n!}\sum_{k=0}^{n}%
\frac{(-n)_k \, x^k}{(\gamma +1 )_k(\alpha+1 )_k} {}_3F_2\left( 
\begin{array}{c}
k-n,\,1-\alpha+k ,\, \gamma\  \\ 
-\alpha-n,\, \gamma+k+1%
\end{array}
\big| 1\right).
\end{equation}
We now give a generating function for the ALPs. It was already found in \cite[p.25]{AW1984} by using the method of differential equations. The alternative proof we present here (see Sect.\ref{proofs} below) is direct and simpler.
\begin{proposition}\label{prop_GF_ALP}
A generating function of the ALPs is 
\begin{eqnarray}\label{gen_ALP_1}
\sum_{n=0}^{+\infty}t^n\mathscr{L}_n^{(\alpha)}(x;\gamma )=(1-t)^{-\gamma-\alpha-1}\exp\left(\frac{x\,t}{t-1}\right)\Phi_1 \left[\gamma,\,\gamma+\alpha,\, \gamma+1;\, \frac{t}{t-1};\, \frac{-x\,t}{t-1} \right]
\end{eqnarray}
where $|t|<\frac{1}{2} ,\, x\in\mathbb{R}$ and $\Phi_1$ is the Humbert confluent hypergeometric function defined by \eqref{Phi_1F1}.\\ In particular, for $\gamma=0$, we recover the generating function of classical Laguerre polynomials
\begin{eqnarray}\label{gener_LP}
\sum_{n=0}^{+\infty}t^n L_n^{(\alpha)}(x)=(1-t)^{-\alpha-1}\exp\left( \frac{x\,t}{t-1} \right),\quad |t|<1 ,\, |x|<+\infty .
\end{eqnarray}
\end{proposition}
Substituting the connection formula (\cite[p.24(2.20)]{AW1984}) $\mathscr{L}_n^{(\alpha)}(x;\gamma )=\sum_{k=0}^{n}\frac{\gamma}{k+\gamma}L_{n-k}^{(\alpha)}(x)L_k^{(-\alpha)}(-x)$ in \eqref{gen_ALP_1}, applying \eqref{RF} to the obtained expression, and then identifying the generating function \eqref{gener_LP} in the final expression to get the following result.
\begin{corollary}
A generating function for Laguerre polynomials
\begin{eqnarray}\label{Cor_gener_ALP_1}
\sum_{n=0}^{+\infty} \frac{\gamma}{n+\gamma}L_n^{(\alpha)}(x)t^n=(1-t)^{-\gamma}\Phi_1 \left[\gamma,\,\gamma -\alpha,\, \gamma+1;\, \frac{t}{t-1};\, \frac{x\,t}{t-1} \right]\quad (\gamma \text{ arbitrary}).
\end{eqnarray}
\end{corollary}
The above generating function generalizes the one given in \cite[Eq(9.12.12)]{askey-scheme}, which we obtain by setting $\gamma=\alpha$ in \eqref{Cor_gener_ALP_1}.
\begin{remark}
Notice that, by a different method, the explicit polynomial form \eqref{ALP_exp}
was found by Askey and Wimp in \cite[p.22, Eq(2.8)]{AW1984} and it can also be obtained from the associated Meixner-Pollaczek
polynomials by using the limit relation: 
\begin{equation*}
\mathscr{L}_n^{(\alpha)}(x;\gamma )=\lim\limits_{\varphi\to 0}\mathscr{P}%
_n^{((\alpha+1)/2)}\left(\frac{-x}{2\sin\varphi};\,\varphi,\, \gamma \right),
\end{equation*}
as was explained by Rahman in \cite[p.7]{Rahman2000}.
\end{remark}
\section{The associated Krawtchouk polynomials}\label{sect_AKP}
The associated Krawtchouk polynomials (AKPs) can be defined by the three-term recurrence relation:
\begin{eqnarray}\label{AKP_RR}
\begin{split}
p(N-n-\gamma)\mathscr{K}_{n+1}(x;p,N,\gamma)=\left[pN+(n+\gamma)(1-2p)-x \right]\mathscr{K}_{n}(x;p,N,\gamma)\\
-(n+\gamma)(1-p)\mathscr{K}_{n-1}(x;p,N,\gamma),
\end{split}
\end{eqnarray}
with initial conditions $\mathscr{K}_{-1}(x;p,N,\gamma)=0$ and $\mathscr{K}_{0}(x;p,N,\gamma)=1$. Applying the criterion \eqref{coef_cond} to \eqref{AKP_RR} we see that the orthogonality is obtained in the following cases:
\begin{enumerate}
\item[(i)] $p<0$, $\quad N-\gamma <0$
\item[(ii)] $0<p<1$, $\quad n=0,1,\cdots, \lfloor N-\gamma \rfloor$
\item[(iii)] $p>1$, $\quad N-\gamma <0$.
\end{enumerate}
The notation $\lfloor x \rfloor$ stands for the largest integer less than or equal to $x$. The AKPs are related to the AMPs in the following way:
\begin{equation}\label{AKP-AMP}
\mathscr{K}_{n}(x;p,N,\gamma)=\frac{\mathscr{M}_{n}(x;-N,p/(p-1),\gamma)}{(-N+\gamma)_n}.
\end{equation}
From \eqref{AMP_4F3} and \eqref{AKP-AMP} follows the explicit representation for the AKPs 
\begin{eqnarray}\label{AKP_4F3_1}
\begin{split}
\mathscr{K}_n(x;p , N, \gamma )=\left(\frac{p-1}{p}\right)^{n}\frac{(\gamma +1 )_n}{n!}\sum_{k=0}^{n}(1-p)^{-k}\frac{(-n)_k(\gamma-N+x)_k}{(\gamma +1 )_k(\gamma-N)_k}\\ 
\times {}_4F_3\left( 
\begin{array}{c}
k-n,\,\gamma-N+x+k ,\, \gamma-N -1,\, \gamma\ \\ 
 \gamma-N +x,\, \gamma-N +k, \, \gamma+1+k
\end{array}
\big| 1\right)
\end{split}
\end{eqnarray} 
valid under restrictions $(i),\ (ii), (iii)$ given above. An alternative form of $\mathscr{K}_n(x;p , N, \gamma )$ stems from \eqref{AMP_4F3_2} and is given by
\begin{eqnarray}\label{AKP_4F3_2}
\begin{split}
\mathscr{K}_n(x;p , N, \gamma )=\frac{(\gamma +1 )_n}{n!}\sum_{k=0}^{n}p^{-k}\frac{(-n)_k(\gamma-x)_k}{(\gamma +1 )_k(\gamma-N)_k} {}_4F_3\left( 
\begin{array}{c}
k-n,\,\gamma-x+k ,\, \gamma-N -1,\, \gamma\ \\ 
 \gamma-x,\, \gamma-N +k, \, \gamma+1+k
\end{array}
\big| 1\right).
\end{split}
\end{eqnarray}
Next, we establish the following generating function for the AKPs.
\begin{theorem}
For $x=x_j$, $j=\{0,1,\cdots , M\}$ (where it is possible that $M=\infty$), we have 
\begin{eqnarray}\label{gen_AKP}
\sum_{n=0}^{+\infty}t^n\mathscr{K}_n(x;p , N, \gamma )=\left(1-t\right)^{-1}F_1\left[1 , \gamma , -x ; \gamma -N ; \frac{t(p-1)}{p}, \frac{t}{p(t-1)} \right],
\end{eqnarray} 
valid for $\max \left\{\left|\frac{t(p-1)}{p}\right|, \left|\frac{t}{p(t-1)}\right| \right\}<1$.
\end{theorem}
It should be noted here that the case when $n=N-\gamma\in\mathbb{N}$ must be understood by continuity. In fact, in this case the hypergeometric ${}_4F_3(1)$ in the expression of $\mathscr{K}_n(x;p , N, \gamma)$ reduces to ${}_3F_2(1)$ and we still have a polynomial of degree $N-\gamma$. For instance, if we take $n=N-\gamma$ in \eqref{AKP_4F3_2} we get
\begin{eqnarray}\label{AKP_4F3_3}
\begin{split}
\mathscr{K}_{N-\gamma}(x;p , N, \gamma )=\frac{(\gamma +1 )_{N-\gamma}}{(N-\gamma)!}\sum_{k=0}^{N-\gamma}p^{-k}\frac{(\gamma-x)_k}{(\gamma +1 )_k} {}_3F_2\left( 
\begin{array}{c}
\gamma-x+k ,\, \gamma-N -1,\, \gamma\ \\ 
 \gamma-x,\, \gamma+1+k
\end{array}
\big| 1\right).
\end{split}
\end{eqnarray}

Therefore, in the case when $N-\gamma\in\mathbb{N}$, we will need a special notation for the generating function of the AKPs. We define the $N$-th partial sum of a power series in $t$ by 
\begin{equation}
\left[f(t)\right]_N:=\sum\limits_{k=0}^{N}\frac{f^{(k)}(0)}{k!}t^k,
\end{equation}
for every function $f$ for which $f^{(k)}(0)$, $k=0,1,\cdots , N$ exists. Then, for $x=0,1,2,\cdots , N-\gamma$, we have the following generating function for the AKPs
\begin{eqnarray}\label{gen_AKP_2}
\sum_{n=0}^{N-\gamma}t^n\mathscr{K}_n(x;p , N, \gamma )=\left[\left(1-t\right)^{-1}F_1\left[1 , \gamma , -x ; \gamma -N ; \frac{t(p-1)}{p}, \frac{t}{p(t-1)} \right]\right]_{N-\gamma}.
\end{eqnarray} 
Taking $\gamma=0$ in the above relation gives the generating function of the classical Krawtchouk polynomials, $K_n(x;p,N)=\mathscr{K}_n(x;p , N, 0)$, obtained for $x=0,1,\cdots , N$ in \cite[(1.10.13)]{askey-scheme}:
\begin{eqnarray}
\sum_{n=0}^{N}t^nK_n(x;p,N)=\left[\left(1-t\right)^{-1}{}_3F_2\left( 
\begin{array}{c}
1 ,\, -x\ \\ 
 -N
\end{array}
\big| \frac{t}{p(t-1)} \right)\right]_{N},
\end{eqnarray} 
where $N$ is a nonnegative integer. 
\section{Concluding remark}
In this paper, we have established explicit representation as well as generating function of the AMPs. These results are then used to derive similar relations for a sequence of associated orthogonal polynomials which belong to the same family by limiting procedures or connection formula. We recall that orthogonality measures of AMPs, ACPs, and AKPs are unknown. An attempt to get the orthogonality measures for the AMPs and ACPs was made in \cite{ILV_1988} where the authors compute the Stieljes transform of the orthogonality measures of these two polynomials. They showed that these measures are unique. As is well known, if elements of one set of orthogonal polynomials converge to those of another set, and the measures are uniquely determined, then there must be the same corresponding limit relation for the measures. This means that it suffices to find the orthogonality measure of the AMPs to deduce those of the other related polynomials. This question may be the subject of forthcoming works.

\section{Proofs}\label{proofs}

This section is devoted to some technical proofs.

\begin{proof}[\textbf{Proof of corollary }\ref{corollary_sum_4F3}]
By combining \eqref{m_n} and \eqref{AMP_4F3_2} then making the modifications $(t,a,b,y)\leftarrow (\tilde{c}, \gamma , \gamma +\beta -1,-x)$, we can readily obtain the following finite sum formula for ${}_4F_3(1)$:
\begin{eqnarray}
\begin{split}
&\frac{(a+1 )_n(b+1)_n}{n!}\sum_{k=0}^{n}t^k\frac{(-n)_k(a+y)_k}{(a +1 )_k(b+1)_k} {}_4F_3\left( 
\begin{array}{c}
k-n,\,a+y+k ,\, a,\, b\ \\ 
 a+y,\, b+1 +k, \, a+1+k
\end{array}
\big| 1\right)\\
&=\frac{1}{b-a}\left\{ (b)_{n+1} \, {}_2F_1\left( 
\begin{array}{c}
1-y,a \\ 
a-b+1
\end{array}
\big|t\right){}_2F_1\left( 
\begin{array}{c}
y,-n-a \\ 
b-a+1
\end{array}
\big|t\right)\right.\\
&\left. -(a)_{n+1}\,{}_2F_1\left( 
\begin{array}{c}
b-a+1-y,b \\ 
b-a+1
\end{array}
\big|t\right){}_2F_1\left( 
\begin{array}{c}
a-b+y,-n-b \\ 
a-b+1
\end{array}
\big|t\right)\right\}.
\end{split}
\end{eqnarray}
Next, to ${}_2F_1$ in the second line of the RHS of the above equation, we apply the Euler transformation \cite[p.33(21)]{Sr-Ma1984}:
\begin{equation}\label{Euler_trans}
{}_2F{}_1\left( 
\begin{array}{c}
A,B \\ 
C
\end{array}
\big|z\right)=(1-z)^{C-A-B}{}_2F_1\left( 
\begin{array}{c}
C-A,C-B \\ 
C
\end{array}
\big|z\right),
\end{equation}
where $C\neq 0,-1,-2,... , |\text{arg}(1-z)|<\pi$, for $A=b-a+1-y, \ B=b,\ C=b-a+1, z=t$ in the first ${}_2F_1$ and $A=a-b+y, \ B=-n-b,\ C=a-b+1, z=t$ to get the formula in \eqref{sum_fin_4F3}. 
\end{proof}
\begin{proof}[\textbf{Proof of Theorem} \ref{thm_GF_AMP}]
Denoting the left-hand side of \eqref{gen2_AMP} by $\Upsilon(x,t)$ and substituting the expression of $\mathscr{M}_n(x;\beta , c, \gamma )$, we obtain
\begin{eqnarray*}
\begin{split}
\Upsilon(x,t)&=\sum_{n=0}^{+\infty}\frac{t^n(\gamma+1)_n}{n!}\sum_{k=0}^{n}(1-c)^k\frac{(-n)_k(\gamma+\beta+x)_k}{(\gamma +1 )_k(\gamma+\beta)_k} {}_4F_3\left( 
\begin{array}{c}
k-n,\,\gamma+\beta+x+k ,\, \gamma+\beta -1,\, \gamma\ \\ 
 \gamma+\beta +x,\, \gamma+\beta +k, \, \gamma+1+k
\end{array}
\big| 1\right)\\
&= \sum_{n=0}^{+\infty}\sum_{k=0}^{n}\sum_{j=0}^{n-k}\frac{t^n(1-c)^k}{n!j!}\frac{(-n)_k(-n-k)_j(\gamma+1)_n(\gamma+\beta+x)_k(\gamma+\beta+x+k)_j (\gamma+\beta -1)_j (\gamma)_j}{(\gamma +1 )_k(\gamma+\beta)_k(\gamma+\beta +x)_j(\gamma+\beta +k)_j(\gamma+1+k)_j}\\
&= \sum_{n=0}^{+\infty}\sum_{k=0}^{+\infty}\sum_{j=0}^{+\infty}\frac{t^{n} (t(c-1))^{k}(-t)^j}{n!j!}\frac{(\gamma+1+k+j)_{n}(\gamma+\beta+x)_k(\gamma+\beta+x+k)_j (\gamma+\beta -1)_j (\gamma)_j}{(\gamma +\beta)_k(\gamma+\beta+x)_j(\gamma+\beta+k)_j },
\end{split}
\end{eqnarray*} 
where we have used, respectively, the series transformation 
\cite[p.102(17)]{Sr-Ma1984}:
\begin{equation}\label{series_trans_3}
\sum_{n=0}^\infty\sum_{k=0}^n \sum\limits_{j=0}^{n-k}A(j,k,n)=\sum_{n=0}^\infty\sum_{k=0}^\infty \sum\limits_{j=0}^{\infty}A(j,k,n+k+j),
\end{equation}
and the identities
\begin{equation}\label{1/n!}
\frac{(-n-k-j)_k(-n-j)_j}{(n+k+j)!}=\frac{(-1)^{k+j}}{n!},
\end{equation}
\begin{equation}
\frac{(\gamma+1)_{n+k+j}}{(\gamma+1)_k(\gamma+1+k)_j}=(\gamma+1+k+j)_{n},
\end{equation}
obtained by virtue of relations
\begin{eqnarray}
\label{Pocha_neg_integ}
&&(-n)_k=(-1)^k\frac{n!}{(n-k)!} ,\quad k=0\leq k \leq n,\\ \label{Pocha_alter}
&&(a)_{n+m}=(a)_n(a+n)_m=(a)_m(a+m)_n, \ a\in\mathbb{C},\ n,m\in\mathbb{N}.
\end{eqnarray}
After summation over $n$ in the last expression of $\Upsilon(x,t)$, we obtain
\begin{eqnarray*}
\Upsilon(x,t)= (1-t)^{-\gamma-1}\sum_{k=0}^{+\infty}\sum_{j=0}^{+\infty}\left(\frac{t(c-1)}{1-t}\right)^{k}\frac{\left(\frac{t}{t-1}\right)^j}{j!}\frac{(\gamma+\beta+x)_k(\gamma+\beta+x+k)_j (\gamma+\beta -1)_j (\gamma)_j}{(\gamma +\beta)_k(\gamma+\beta+x)_j(\gamma+\beta+k)_j }.
\end{eqnarray*} 
We apply again the identity \eqref{Pocha_alter} to get, after simplifications, 
\begin{eqnarray*}
\Upsilon(x,t)&=& (1-t)^{-\gamma-1}\sum_{k=0}^{+\infty}\sum_{j=0}^{+\infty}\left(\frac{t(c-1)}{1-t}\right)^{k}\frac{\left(\frac{t}{t-1}\right)^j}{j!}\frac{(\gamma+\beta+x+j)_k (\gamma+\beta -1)_j (\gamma)_j}{(\gamma +\beta)_j(\gamma+\beta+j)_k}\\
&=&(1-t)^{-\gamma-1}\sum_{j=0}^{+\infty}\frac{(\gamma+\beta-1)_j(\gamma)_j}{(\gamma+\beta)_j}  {}_2F_1\left( 
\begin{array}{c}
\gamma +\beta+x+ j,\, 1\ \\ 
 \gamma+\beta+j
\end{array}
\left| \frac{t(c-1)}{1-t}\right.\right)\frac{\left(\frac{t}{t-1}\right)^j}{j!}.
\end{eqnarray*} 
Next, to ${}_2F_1$ in the above equation, we apply the Euler transformation \cite[p.33(21)]{Sr-Ma1984}:
\begin{equation}\label{Euler_trans}
{}_2F{}_1\left( 
\begin{array}{c}
a,b \\ 
e
\end{array}
\big|z\right)=(1-z)^{e-a-b}{}_2F_1\left( 
\begin{array}{c}
e-a,e-b \\ 
e
\end{array}
\big|z\right),
\end{equation}
where $e\neq 0,-1,-2,... , |\text{arg}(1-z)|<\pi$, to get
\begin{eqnarray*}
\Upsilon(x,t)=(1-t)^{x-\gamma}(1-ct)^{-x-1}\sum_{j=0}^{+\infty}\frac{(\gamma+\beta-1)_j(\gamma)_j}{(\gamma+\beta)_j}  {}_2F_1\left( 
\begin{array}{c}
-x,\, \gamma +\beta-1+ j \ \\ 
 \gamma+\beta+j
\end{array}
\left| \frac{t(c-1)}{1-t}\right.\right)\frac{\left(\frac{t}{t-1}\right)^j}{j!}.
\end{eqnarray*}
We identify the infinite series as the Appell hypergeometric function $F_1$ defined by (\cite[p.53(4)]{Sr-Ma1984}):
\begin{eqnarray}\label{Appel_F1}
\begin{split}
F_1\left[ \alpha ,\, \beta_1 ,\, \beta_2 ;\, \sigma ;\, x,\,y \right]&=\sum\limits_{m,n=0}^{+\infty} \frac{(\alpha)_{m+n}(\beta_1)_m(\beta_2)_n}{(\sigma)_{m+n}}\frac{x^m}{m!}\frac{y^n}{n!}\\
&=\sum\limits_{m=0}^{+\infty} \frac{(\alpha)_{m}(\beta_1)_m}{(\sigma)_{m}} {}_2F_1\left( 
\begin{array}{c}
\alpha+m , \beta_2 \\ 
\sigma+m
\end{array}
\big|y\right)\frac{x^m}{m!}, \quad \max\{|x|,|y|\}<1.
\end{split}
\end{eqnarray}
Thus, we have
\begin{eqnarray*}
\Upsilon(x,t)=(1-t)^{x-\gamma}(1-ct)^{-x-1}F_1\left[\gamma+\beta-1 ; \gamma , -x ; \gamma +\beta ; \frac{t}{t-1}, \frac{t(1-c)}{t-1} \right].
\end{eqnarray*}
To get \eqref{gen2_AMP}, it suffices to apply the transformation (see \cite[p.78]{Bailey1964}):
\begin{eqnarray}
F_1\left[ \alpha ,\, \beta_1 ,\, \beta_2 ;\, \sigma ;\, x,\,y \right]=(1-x)^{-\beta_1}(1-y)^{-\beta_2}F_1\left[ \sigma -\alpha ,\, \beta_1 ,\, \beta_2 ;\, \sigma ;\, \frac{x}{x-1},\,\frac{y}{y-1} \right]\label{F1_trans1}
\end{eqnarray}
Finally, the formula \eqref{gener2_MPP} is a particular case of \cite[Eq(9.10.13)]{askey-scheme} and it is obtained here directly by putting $\gamma=0$ in \eqref{gen2_AMP}. The proof is complete. 
\end{proof}
%%%%%%%%%%%%%%%%%%%%%%%%%%%%%%%%%%%%%%%%%%%%%%%%%%%%%%%%%%%%%%%%%%%%%%%%%%%%%%%%%%%%%%%%%%%
%\newpage
%%%%%%%%%%%%%%%%%%%%%%%%%%%%%%%%%%%%%%%%%%%%%%%%%%%%%%%%%%%%%%%%%%%%%%%%%%%%%%%%%%%%%%%%%%%%%
\begin{proof}[\textbf{Proof of Proposition} \ref{prop_GF_ACP}]
In fact, recalling the expression of \eqref{ACP_exp} and denoting the left-hand side of \eqref{gener_ACP} by $\mathscr{G}(x,t)$, we obtain
\begin{eqnarray}
\mathscr{G}(x,t)&=&\sum_{n=0}^{+\infty}\frac{t^n}{n!}\sum_{k=0}^{n}(-a)^{-k}\frac{(-n)_k(\gamma-x)_k}{(\gamma +1 )_k}
{}_3F_2\left( 
\begin{array}{c}
k-n,\,\gamma-x+k ,\, \gamma\ \\ 
 \gamma-x,\, \gamma+k+1
\end{array}
\big| 1\right)\\
&=& \sum_{n=0}^{+\infty}\sum_{k=0}^{n}\sum_{j=0}^{n-k}\frac{t^n (-a)^{-k}}{n!j!}\frac{(-n)_k(k-n)_j(\gamma-x)_k(\gamma-x+k)_j(\gamma)_j}{(\gamma +1 )_k(\gamma-x)_j(\gamma+k+1)_j }\\
&=& \sum_{n=0}^{+\infty}\sum_{k=0}^{+\infty}\sum_{j=0}^{+\infty}\frac{t^{n+k+j} (-a)^{-k}}{(n+k+j)!j!}\frac{(-n-k-j)_k(-n-j)_j(\gamma-x)_k(\gamma-x+k)_j(\gamma)_j}{(\gamma +1 )_k(\gamma-x)_j(\gamma+k+1)_j }.
\end{eqnarray} 
These calculations can be done by applying the identity \eqref{series_trans_3}. Next, applying \eqref{1/n!} and the relation
$(\gamma +y)_k(\gamma +y+k)_j=(\gamma +y)_j(\gamma +y+j)_k$, for $y=-x,\ 1$, to the last expression of $\mathscr{G}(x,t)$ to get
\begin{eqnarray}
\mathscr{G}(x,t)&=& e^t\sum_{k=0}^{+\infty}\sum_{j=0}^{+\infty}\frac{(\gamma-x+j)_k(\gamma)_j}{(\gamma+1+j)_k(\gamma +1 )_j }\left(\frac{t}{a}\right)^k\frac{  (-t)^j}{j!}  \nonumber \\ \label{2F1_G(x,t)}
&=& e^t\sum_{j=0}^{+\infty}\frac{(\gamma)_j}{(\gamma +1 )_j } {}_2F_1\left( 
\begin{array}{c}
\gamma-x+j ,\, 1\ \\ 
 \gamma+1+j
\end{array}
\left| \frac{t}{a}\right.\right)\frac{  (-t)^j}{j!}.
\end{eqnarray}
Next, to ${}_2F_1$ in \eqref{2F1_G(x,t)} we apply the Euler transformation \eqref{Euler_trans}, to get
\begin{eqnarray}\label{G_sum_2F1}
\mathscr{G}(x,t)=e^t\left(1-\frac{t}{a}\right)^x\sum_{j=0}^{+\infty}\frac{(\gamma)_j}{(\gamma +1 )_j } {}_2F_1\left( 
\begin{array}{c}
\gamma +j,\, x+1\ \\ 
 \gamma+1+j
\end{array}
\left| \frac{t}{a}\right.\right)\frac{(-t)^j}{j!}.
\end{eqnarray}
By recognizing the Humbert confluent hypergeometric function $\Phi_1$ defined by \cite[p.58(36)]{Sr-Ma1984}:
\begin{eqnarray}\label{Phi_1F1}
\begin{split}
\Phi_1\left[\alpha_1 , \,\lambda ;\,\alpha_2 ;\, x,\, y \right]&=\sum\limits_{m,n=0}^{+\infty} \frac{(\alpha_1)_{m+n}(\lambda)_m}{(\alpha_2)_{m+n}}\frac{x^m}{m!}\frac{y^n}{n!}\\ 
&=\sum\limits_{n=0}^{+\infty} \frac{(\alpha_1)_{n}}{(\alpha_2)_{n}} {}_2F_1\left( 
\begin{array}{c}
\alpha_1+n , \ \lambda \\ 
\alpha_2+n
\end{array}
\big|x\right)\frac{y^n}{n!}, \quad |x|<1,\, |y|<+\infty .
\end{split}
\end{eqnarray}
in the right-hand side of \eqref{G_sum_2F1}, we complete the proof of the Proposition \ref{prop_GF_ACP}.
\end{proof}
%%%%%%%%%%%%%%%%%%%%%%%%%%%%%%%%%%%%%%%%%%%%%%%%%%%%%%%%%%%%%%%%%%%%%%%%%%%%%%%%%%%%
\begin{proof}[\textbf{Proof of Proposition} \ref{prop_GF_ALP}]
Denote the left-hand side of \eqref{gen_ALP_1} by $\Lambda(x,t)$. Similar calculations to the ones of the proof of Proposition \ref{prop_GF_ACP} give
\begin{eqnarray}
\Lambda(x,t)&=&\sum_{n=0}^{+\infty}\frac{t^n}{n!}\sum_{k=0}^{n}\frac{(\gamma +\alpha+1 )_n(-n)_k \, x^k}{(\gamma +1 )_k(\gamma +\alpha+1 )_k}
{}_3F_2\left( 
\begin{array}{c}
k-n,\,\gamma+\alpha ,\, \gamma\ \\ 
 \gamma+\alpha+k+1,\, \gamma+k+1
\end{array}
\big| 1\right)\\
%&=& \sum_{n=0}^{+\infty}\sum_{k=0}^{n}\sum_{j=0}^{n-k}\frac{t^n x^{k}}{n!j!}\frac{(\gamma +\alpha+1 )_n(-n)_k(k-n)_j(\gamma+\alpha)_j(\gamma)_j}{(\gamma +1 )_k(\gamma +\alpha+1 )_k(\gamma+\alpha+k+1)_j(\gamma+k+1)_j }\\
%&=& \sum_{n=0}^{+\infty}\sum_{k=0}^{+\infty}\sum_{j=0}^{+\infty} \frac{t^{n+k+j} x^{k}}{(n+k+j)!j!}\frac{(\gamma +\alpha+1 )_{n+k+j}(-n-k-j)_k(-n-j)_j(\gamma+\alpha)_j(\gamma)_j}{(\gamma +1 )_k(\gamma +\alpha+1 )_k(\gamma+\alpha+k+1)_j(\gamma+k+1)_j }\\
&=& \sum_{n=0}^{+\infty}\sum_{k=0}^{+\infty}\sum_{j=0}^{+\infty} \frac{(\gamma +\alpha+1 )_{n+k+j}t^{n} }{n!}\frac{(-xt)^{k}(-t)^j(\gamma+\alpha)_j(\gamma)_j}{j!(\gamma +1 )_{k+j}(\gamma +\alpha+1 )_{k+j} }\\
&=& \sum_{n=0}^{+\infty}\sum_{k=0}^{+\infty}\sum_{j=0}^{+\infty} \frac{(\gamma +\alpha+1 +k+j)_{n}t^{n} }{n!}\frac{(-xt)^{k}(-t)^j(\gamma+\alpha)_j(\gamma)_j}{j!(\gamma +1 )_{k+j} }\\
&=& (1-t)^{-\gamma-\alpha -1} \sum_{k=0}^{+\infty}\sum_{j=0}^{+\infty} \frac{(\frac{xt}{t-1})^{k}(\frac{t}{t-1})^j(\gamma+\alpha)_j(\gamma)_j}{j!(\gamma +1 )_{k+j} }\\
&=& (1-t)^{-\gamma-\alpha -1} \sum_{j=0}^{+\infty} \frac{(\frac{t}{t-1})^j(\gamma+\alpha)_j(\gamma)_j}{j!(\gamma +1 )_{j} } {}_1F_1\left( 
\begin{array}{c}
1\ \\ 
 \gamma+1+j
\end{array}
\left| \frac{xt}{t-1}\right.\right).
\end{eqnarray} 
Next, we use the Kummer transformation \cite[p.37(7)]{Sr-Ma1984}:
\begin{equation}\label{Kummer_trans}
{}_1F_1\left( a; b; z\right)=e^z\, {}_1F_1\left( b-a; b; -z\right),
\end{equation}
to get
\begin{eqnarray}
\Lambda(x,t)=(1-t)^{-\gamma-\alpha -1}\exp\left(\frac{xt}{t-1}\right) \sum_{j=0}^{+\infty} \frac{(\frac{t}{t-1})^j(\gamma+\alpha)_j(\gamma)_j}{j!(\gamma +1 )_{j} } {}_1F_1\left( 
\begin{array}{c}
\gamma +j\ \\ 
 \gamma+1+j
\end{array}
\left| \frac{-xt}{t-1}\right.\right).
\end{eqnarray}
Identification of the summation over $j$ as the Humbert confluent hypergeometric function $\Phi_1$ then completes the proof of Proposition \ref{prop_GF_ALP}.
\end{proof}

%%%%%%%%%%%%%%%%%%%%%%%%%%%%%%%%%%%%%%%%%%%%%%%%%%%%%%%%%%%%%%%%%%%%%%%%%%%%%%%%%%%%%%%%%%
\subsection*{Acknowledgments} I would like to thank Professor Zouha\"ir Mouayn for his comments and invaluable suggestions that improved the presentation of this manuscript. I would like also to thank the Moroccan Association of Harmonic Analysis and Spectral Geometry.    
\subsection*{Data availability statement} All data generated or analyzed during this study are included in this published article.

\end{document}